\documentclass[twoside]{article}
\newcommand{\mytitle}{Equidistribution of Algebraic Numbers of Norm One
  in Quadratic Number Fields} 
\newcommand{\keywords}{Weil height, equidistribution, Hecke
  $L$-function, visible points, quadratic algebraic numbers}

\usepackage{amsmath,amssymb, amsthm,amsfonts,amscd,mathrsfs,pifont} 
\usepackage{eucal}  
\usepackage{verbatim}     
\usepackage[dvips]{graphicx}
\usepackage{hyperref}

\newtheorem{thm}{Theorem}[section]

\newtheorem{lemma}[thm]{Lemma}

\newtheorem{thm*}{Theorem}[]
\newtheorem{cor*}[thm*]{Corollary}
\newtheorem{claim*}[thm*]{Claim}
\newtheorem{lemma*}[thm*]{Lemma}
\newtheorem{prop*}[thm*]{Proposition}
\newtheorem{conj*}[thm*]{Conjecture}

\theoremstyle{definition}

\newtheorem{question*}{Question}[section]

\newtheorem{defn*}{Definition}

\theoremstyle{remark}

\newcommand{\BB}[1]{\ensuremath{\mathbb{#1}}}

\newcommand{\N}{\ensuremath{\BB{N}}}
\newcommand{\Q}{\ensuremath{\BB{Q}}}
\newcommand{\R}{\ensuremath{\BB{R}}}
\newcommand{\Z}{\ensuremath{\BB{Z}}}
\newcommand{\C}{\ensuremath{\BB{C}}}
\newcommand{\T}{\ensuremath{\BB{T}}}

\newcommand{\Ox}{\ensuremath{\mathcal{O}}}

\newcommand{\mf}{\ensuremath{\mathfrak}}

\newcommand{\wt}{\ensuremath{\widetilde}}

\newcommand{\ip}[1]{\mathrm{Im}(#1)}
\newcommand{\rp}[1]{\mathrm{Re}(#1)}

\DeclareMathOperator{\Res}{Res}

\usepackage{ifpdf}
\ifpdf
\fi

\usepackage[usenames]{color}

\def\Ox{{\mathcal O}}

\newcounter{nootje}
\numberwithin{equation}{section} 

\numberwithin{equation}{section}

\begin{document}
\title{\bf \mytitle}  
\author{\sc Kathleen L.~Petersen \rm{and} \sc Christopher
  D.~Sinclair\footnote{The second author was supported in part by the
    National Science Foundation (DMS-0801243)} } 
\maketitle

\begin{abstract}
  Given a fixed quadratic extension $K$ of $\Q$, we consider the
  distribution of elements in $K$ of norm 1 (denoted $\mathcal
  N$).  When $K$ is an imaginary quadratic extension, $\mathcal N$ is
  naturally embedded in the unit circle in $\C$ and we show that it is
  equidistributed with respect to inclusion as ordered by the
  absolute Weil height.  By Hilbert's Theorem 90, an element in
  $\mathcal N$ can be written as $\alpha/\overline{\alpha}$ for some
  $\alpha \in \Ox_K$, which yields another ordering of
  $\mathcal N$ given by the minimal norm of the associated algebraic
  integers.  When $K$ is imaginary we also show that $\mathcal N$ is
  equidistributed in the unit circle under this norm ordering.  When $K$ is a
  real quadratic extension, we show that $\mathcal N$ is
  equidistributed with respect to norm, under the map $\beta \mapsto
  \log| \beta | \bmod{\log | \epsilon^2 |}$ where $\epsilon$ is a
  fundamental unit of $\mathcal O_K$.
\end{abstract}


\section{Introduction}

Let $d$ be a square-free integer and let $K = \Q(\sqrt{d})$.  We may
embed $K$ in $\C$ or $\R$ (depending on whether $d$ is negative or
positive, respectively) by setting
\[
K = \{ a + b \sqrt{d} : a, b \in \Q \}.
\]
(Where, for concreteness, we take $\sqrt{d} > 0$ if $d > 0$ and
$\ip{\sqrt{d}} > 0$ if $d < 0$.)
We define 
\[ \omega = \left\{ \begin{array}{cl} \sqrt{d} & \mbox{ if }
    d\equiv 2,3 \bmod 4 \\  {\displaystyle \frac{1+\sqrt{d}}{2}} &
    \mbox{ if } d \equiv 1 \bmod 4 \end{array} \right. 
\]  
so that the ring of integers $\mathcal O$ of $K$ is identified with
$\Z \oplus \omega \Z$. 
As usual, if $\beta = a + b \omega$, we define  the conjugate of
$\beta$ by $\overline{\beta} = a + b \overline{\omega}$ and the {\em norm} of $\beta$ by $N(\beta) = N_{K/\Q}(\beta) = \beta
\overline{\beta} = a^2 + ab(\omega + \overline{\omega}) + b^2 \omega
\overline{\omega}$.    

We are interested in the distribution  of the set  
\[
\mathcal N = \{ \beta \in K : N(\beta) = 1 \}.
\]
This is related to the geometry of certain subsets of $\Ox$ 
as embedded in $\R^2$ or $\C$.  When $K$ is an imaginary quadratic
$\mathcal N = K \cap \T$ where $\T \subset \C$ is the unit circle.
When $K$ is real, the map $\lambda:K\rightarrow \R^2$ defined by
$\lambda(\beta)= (\beta,\overline{\beta})$ identifies $\Ox$ with a
lattice in $\R^2$.  The set $\mathcal N$ maps to the intersection of
$\lambda(K)$ with the hyperbola $ \{ (x,y)\in \R^2: xy=1\}$.

The authors would like to thank Ted Chinburg whose questions motivated this work,
and  Ram Murty who supplied the idea behind the proof of Theorem~\ref{thm:4}.

\subsection{The Absolute Weil Height}

Our principle equidistribution result for imaginary quadratics
concerns $\mathcal{N}$ ordered by the absolute Weil
height; we review necessary definitions here.

Two absolute values $| \cdot |_1$ and $| \cdot |_2$ on $K$ (which in
this section can be an arbitrary number field) are equivalent if there
exists some exponent $e > 0$ such that $| \beta |_2 = | \beta |_1^e$
for all $\beta \in K$.  An equivalence class of non-trivial absolute
values is known as a {\em place} of $K$.

When $K = \Q$, we are familiar with many inequivalent absolute
values.  In particular, we set $| \cdot |_{\infty}$ to be the usual
absolute value inherited from $\R$ and, for each rational prime $p$,
we define $| \cdot |_p$ to be the $p$-adic absolute value normalized
so that $| p |_p = p^{-1}$.  A celebrated theorem of Ostrowski shows
that this is a complete set of representatives of the places of $\Q$.
\begin{thm}[Ostrowski]
  Any nontrivial absolute value on $\Q$ is equivalent to exactly one
  $| \cdot |_p$ for $p \in \{\infty, 2, 3, 5, \cdots \}$.
\end{thm}
\noindent (See, for instance, \cite[\S 9.3]{Jacobson:1989kx} for discussion and
a proof of Ostrowski's Theorem).

Each absolute value on $K$ induces an absolute value on $\Q$ by
inclusion.  If $u$ is the place of $\Q$ containing $| \cdot |_p$ and
$v$ is a place of $K$ whose elements when restricted to $\Q$ lie in
$u$, then we will say that $v$ {\em divides} $p$ and write $v | p$. In this
situation, there is a one-to-one correspondence between absolute
values in $u$ and absolute values in $v$.  We remark that there may be
multiple places of $K$ which divide a given $p$.  

It is useful to distinguish a representative in each place of $K$.  If
$v | \infty$, then we will say that $v$ is an {\em Archimedean} place
of $K$ and set $\| \cdot \|_{v} \in v$ to be the unique absolute value
whose restriction to $\Q$ is $| \cdot |_{\infty}$.  If on the other
hand, $p$ is a (finite) rational prime and $v | p$ then we say $v$ is
{\em non-Archimedean} and write $\| \cdot \|_v \in v$ for the
absolute value which restricts to the usual $p$-adic absolute value on
$\Q$.

We write $K_v $ for the completion of $K$ with respect to any absolute 
value in $v$, and define the {\em local degree} of $K$ with respect to
$v$ by $n_v = [K_v : \Q_p ]$. It can be shown that
\[
\sum_{v | p} n_v = [K : \Q] = n,
\]
where $n = [K : \Q]$.  

For each place $v$ of $K$ we define the absolute value $| \cdot |_v$
by setting $| \cdot |_v = \| \cdot \|_v^{n_v/n}$.  Using this absolute
value, we define the {\em local height} on $K_v$ by 
\[
H_v(\beta) = \max\{1, | \beta |_v \}, \qquad \beta \in K_v.  
\]
The {\em absolute Weil height} on $K$ is then defined by
\begin{equation}
\label{eq:1}
H(\beta) = \prod_v H_v(\beta) = \prod_{v} \max\{1, | \beta |_v \},
\qquad \beta \in K.
\end{equation}
It is remarkable (though not immediately obvious) that $H(\beta)$ 
is independent of $K$.  That is, if $\beta \in K_1 \cap K_2$ for
distinct number fields $K_1$ and $K_2$ and $H_1$ and $H_2$ are defined
as in (\ref{eq:1}) using the places of $K_1$ and $K_2$ respectively,
then $H_1(\beta) = H_2(\beta)$.  

The set of elements in $K$ of bounded height is finite.  That is, $K$
has the {\em Northcott property} with respect to the absolute Weil
height.  In fact, something much stronger is true:  The set of {\em
  all} algebraic numbers whose degree and height are both less than
some fixed bounds is finite.

\subsection{Equidistribution and Weyl's Criterion}

Perhaps the most classical setting in which to study distribution is
in $\R / \Z$. Suppose $(a_n)_{n=1}^{\infty}$ is a sequence of real numbers.
  Then, $(a_n + \Z)_{n=1}^{\infty}$ is {\em
  equidistributed} in $\R / \Z$ if for all subintervals $[a,b)$ of
$[0,1)$ 
\[ 
\lim_{N\rightarrow \infty} \frac1N  \# \{n\leq N : \{a_n\} \in [a,b)\}  =
b-a
\]
where $\{a\}$ denotes the fractional part of $a$. 
That is, if in the limit, the number of elements in the sequence which
are in the subinterval is dependent only on the length of the
subinterval.  It is worth remarking that equidistribution is dependent
on the ordering of $(a_n)_{n=1}^{\infty}$.  

More generally, suppose $T$ is a compact abelian group with normalized
(probability) 
Haar measure $\mu$.  If $\mathcal H(1), \mathcal H(2), \ldots$ are
finite sets such that $\mathcal H(1) \subseteq \mathcal H(2) \subseteq 
\mathcal H(3) \subseteq \cdots $, then we say $\bigcup_{n=1}^{\infty}
\mathcal H(n)$ is {\em equidistributed} in $T$ if for every continuous function
$f$ on $T$,
\[
\lim_{n \rightarrow \infty} \frac{1}{\#
\mathcal  H(n)} \sum_{b \in \mathcal H(n)} f(b) = \int_T f \, d \mu.
\]

Hermann Weyl gave an equivalent condition formulated in terms of the
characters of $T$.  
\begin{thm}[Weyl]
$\bigcup_{n=1}^{\infty} \mathcal H(n)$ is equidistributed in $T$ if for
every character $\chi$ of T, 
\[
\lim_{n \rightarrow \infty} \frac{1}{\# \mathcal H(n)} \sum_{b \in
  \mathcal H(n)} \chi(b) = \int_T \chi \, d\mu.
\]
\end{thm}
(In fact, Weyl formulated this criterion in \cite{MR1511862} for
equidistribution in $\R / \Z$, which was then extended to
arbitrary compact groups by Eckmann \cite{MR0011302}.)

\section{Statement of Results}

Returning to the case where $K$ is a quadratic extension of $\Q$, for
any $t > 0$ we define 
\[
 \mathcal H(t) = \{ \beta \in \mathcal N : H(\beta) \leq t \}.
\]
Clearly, if $t_1 < t_2$, then $\mathcal H(t_1) \subseteq \mathcal
H(t_2)$, and 
\begin{equation}
\label{eq:7}
\mathcal H = \bigcup_{n=1}^{\infty} \mathcal H(n). 
\end{equation}
(Of course, as a {\em set}, $\mathcal H = \mathcal N$, we use the
notation to distinguish a partial ordering of the underlying set).  
By the Northcott property of $K$,  $\mathcal
H(t)$ is a finite set 
for all $t > 0$.   If $K$ is an imaginary quadratic extension of $\Q$,
we may therefore talk about the equidistribution of
$ \mathcal H$ in $\T$.  

\begin{thm}
\label{thm:1}
Suppose $K$ is an imaginary quadratic extension of $\Q$.  Then,
$\mathcal H$ is equidistributed in $\T$. 
\end{thm}

When $K$ is a real quadratic extension of $\Q$ we have to make some
alterations as $\mathcal N$ is not naturally embedded in a compact
abelian group.  We may reconcile this issue by using Dirichlet's Unit
Theorem.  Clearly $\mathcal N$ is a group under multiplication, and we
define the homomorphism $\Lambda: \mathcal N \rightarrow \R$ by
\[\Lambda(\beta) = \log|\beta|.\]  The proof of Dirichlet's Unit Theorem
shows that the image of the unit group $U$ of $K$ under $\Lambda$ is a
rank 1 lattice, and that $\ker \Lambda = \{1, -1\}$. 
Since $U^2 = \{ \upsilon^2 : \upsilon \in U\}$ is a subgroup of $U$, it follows that
$\Lambda(U^2)$ is likewise a rank 1 lattice in $\R$.  If $\epsilon$ is a 
fundamental unit of $K$, then $U^2$ is generated by $\epsilon^2$ and
$\Lambda(U^2) = \{ n \log|\epsilon^2| : n \in \Z \}$.  Clearly
$\mathbb{S} =  \R \big/ \Lambda(U^2)$ is isomorphic to the 
circle group, and the map \begin{equation}
\label{eq:11}
\varphi(\beta) = \log|\beta| \bmod \log|
\epsilon^2 |
\end{equation}
maps $\mathcal N$ into $\mathbb{S}$.  Given a natural ordering of
$\varphi(\mathcal N)$ we will be in position to appeal to Weyl's
criterion in order to establish an equidistribution result for
elements of norm 1 in real quadratic extensions of $\Q$.  

\subsection{Ordering by Norm}

For any quadratic extension $K$ of $\Q$,
Hilbert's Theorem 90 implies that if $\beta \in \mathcal{N}$, then
there exists an algebraic integer $\alpha \in \mathcal O$ 
such that $\beta = \alpha/\overline{\alpha}$.  That is, the map 
\begin{equation}
\label{eq:12} 
\pi: z \mapsto z/\overline z
\end{equation}
is a homomorphism from $\mathcal O$ onto
$\mathcal N$.  It is easily seen that $\pi^{-1}(1) = \Z$ and
$\pi^{-1}(-1) = \sqrt{d} \Z$.  

We may specify representatives of cosets in $\mathcal O / \ker \pi$ by
restricting to $\alpha = a + b \omega$ where $a, b$ in $\Z$ are
relatively prime and (as an ordered pair) $(a, b)$ is either equal to
$(1,0)$ or satisfies $b > 0$.  We define the set of {\em visible
  points} in $K$ by 
\begin{equation}
\label{eq:9}
\mathcal V = \{ a + b \omega \in \mathcal O :
\mathrm{gcd}(a,b) = 1 \}.
\end{equation}
The `visible' moniker stems from the fact that, the elements of
$\mathcal V$ are exactly those points in $\mathcal O$ which are
`visible' from the origin in the relevant embedding into $\C$ or
$\R^2$.  (That is, an integer is a visible point if the line segment
connecting it to the origin contains no other integers). 
\begin{figure}[h!]
\centering 
\includegraphics[scale=.6]{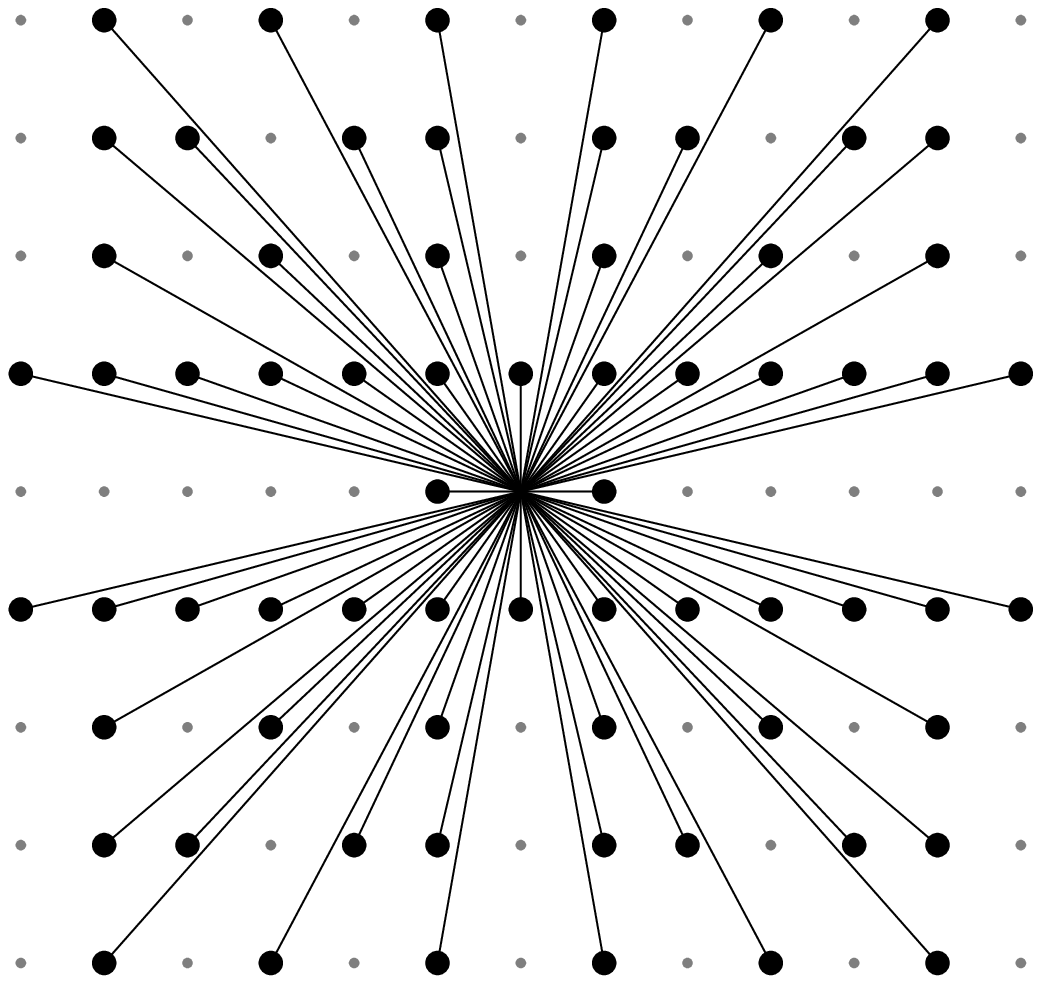}
\includegraphics[scale=.6]{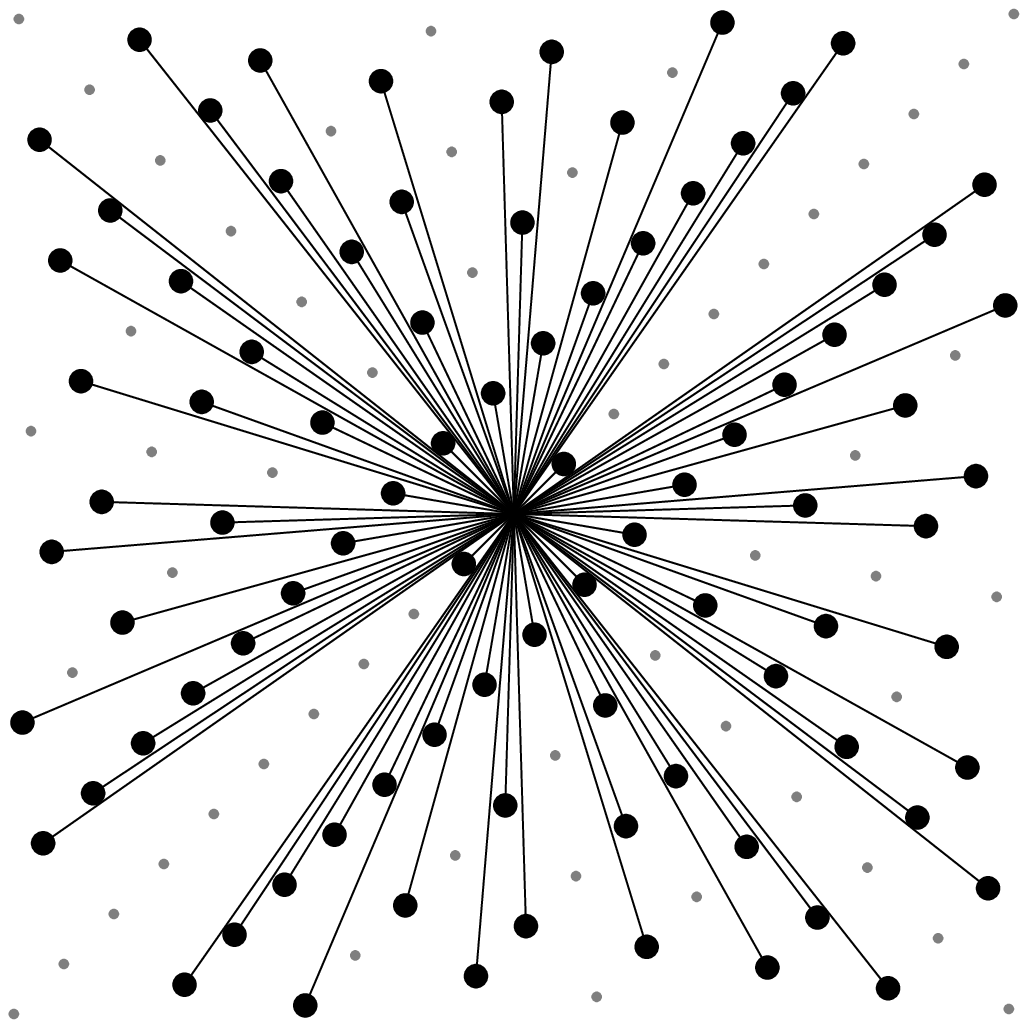}
\begin{caption}{Visible points for $\Q(\sqrt{-2})$ and $\Q(\sqrt{2})$
    as embedded in $\C$ and $\R^2$ respectively.} 
\label{fig:vis}
\end{caption}
\end{figure}

\subsubsection{Imaginary Quadratics}

When $K$ is imaginary, $\pi$ restricted to $\mathcal V$ gives a
two-to-one map onto $\mathcal N$.  Specifically, if $\pi(\alpha) =
\beta$ for some $\alpha \in 
\mathcal V$, then $\pi^{-1}(\beta) = \pm \alpha$.  Since $N(-\alpha) =
N(\alpha)$, we may unambiguously speak of the norm of 
$\pi^{-1}(\beta)$. This gives us another ordering for $\mathcal
N$.  Let
\[
\mathcal M(t) = \{ \beta \in \mathcal N : N\big(\pi^{-1}(\beta)\big) \leq t \},
\]
and
\begin{equation}
\label{eq:8}
\mathcal M = \bigcup_{n=1}^{\infty} \mathcal M(n).
\end{equation}
\begin{thm}
\label{thm:3}
Suppose $K$ be an imaginary quadratic extension of $\Q$.  Then, $\mathcal
M$ is equidistributed in $\T$.   
\end{thm}

\subsubsection{Real Quadratic Extensions}

It is easily seen that the group of units $U$ is contained in
$\mathcal V$.  Moreover we have the following lemma:
\begin{lemma}
\label{lemma:1}
 If $\alpha \in \mathcal V$ and $\upsilon \in U$, then
  $\upsilon \alpha \in \mathcal V$.  
\end{lemma}
\begin{proof}
Case: $d \not\equiv 1 \bmod 4$:
Write $\alpha = a + b\sqrt{d}$ and $\upsilon = u + v \sqrt d
$ so that 
\[
 \upsilon \alpha = (a u + b v d) + (a v + b u) \sqrt{d}.
\]
If $m \neq 0$ is such that 
\[
m \; | \; av + bu \qquad \mbox{and} \qquad m \; | \; au + bv d,
\]
then 
\[
m \; | \; (av + bu) d v \qquad \mbox{and} \qquad m \; | \; (au + bvd) u
\qquad \Rightarrow \qquad m \; | \; a(u^2 - v^2 d)
\]
and 
\[
m \; | \; (av + bu) u \qquad \mbox{and} \qquad m \; | \; (au + bvd) v
\qquad \Rightarrow \qquad m \; | \; b(u^2 - v^2 d).
\]
But, since $u^2 - v^2 d = N(\upsilon) = \pm 1$, and $a$ and $b$ are
relatively prime, $m$ must equal $\pm 1$.  That is, $\upsilon \alpha$
is a visible point.  

Case $d \equiv 1 \bmod 4$:  Write $\alpha = a + b \omega$
and $\upsilon = u + v \omega$, and let $D = (1-d)/4$.
Then $\upsilon \alpha = au + bvD + (av + bu + bv) \omega.$  If there exists $m \neq 0$ such that
\[
m \, | \, a u + b v D \qquad \mbox{and} \qquad m \, | \, av + bu + bv,
\]
then
\[
m \, | \, (au + bv D) v \qquad \mbox{and} \qquad m \, | \, (av + bu + bv) u
\qquad \Rightarrow \qquad m \, | \, b(u^2 + uv + v^2 D)
\]
and
\[
m \, | \, (a u + b v D) (u + v) \qquad \mbox{and} \qquad m \, | \, (a
v + b u + b v) D v \qquad \Rightarrow \qquad m \, | \, a (u^2 + uv +
v^2 D).
\]
But,  $u^2 + uv + v^2 D = N(\upsilon) = \pm 1$, and hence as
before $m = \pm 1$ and $\alpha \in \mathcal V$.  
\end{proof}

Lemma~\ref{lemma:1} holds trivially in the imaginary quadratic
setting, and in both settings we may use it to define an equivalence
relation $\sim$ on $\Ox$ by writing $\gamma \sim \gamma'$ if
there exists a unit $\upsilon \in U$ and positive integers $m$ and
$n\in \N$ such that $\upsilon m \gamma' = n \gamma$.  Geometrically,
$\gamma \sim \gamma'$ if there exists a unit $\upsilon$ such that,
under the appropriate embedding, $\upsilon \gamma$ and $\gamma'$ lie
on the same ray emanating from the origin.  

We denote the set of equivalence classes of by 
$\widetilde{\mathcal O}$.  Notice that with 
\[
\pi(\gamma) = \pi(m \gamma) = \pi( \upsilon m \gamma') = 
\frac{\upsilon^2}{N(\upsilon)} \pi(\gamma') = \pm \upsilon^2 \pi(\gamma')
\]
and hence, with $\varphi$ as in (\ref{eq:11}), 
\[
\varphi \circ \pi(\gamma') = \varphi \circ \pi(\gamma).
\]
That is, we may regard $\varphi \circ \pi$ as a well-defined function on
$\widetilde{\mathcal O}$.   

\begin{lemma}
The function $\varphi \circ \pi$ is two-to-one from $\widetilde{\mathcal
  O}$ into $\mathbb S$.  Specifically, if $\gamma' \not \sim \gamma$
and 
\[
\varphi \circ \pi( \gamma ) = \varphi \circ \pi( \gamma' ),
\]
then $\gamma \sim  \sqrt{d} \gamma'$.  
\end{lemma}
\begin{proof}
Viewing $\varphi(\beta) = \log|\beta| \bmod \log |\epsilon^2|$ as a
homomorphism from $\mathcal N$ to $\mathbb{S}$, it is clear that\[
\ker \varphi =U^2.
\]
Now, if $\gamma$ and $\gamma'$ are such that $\pi(\gamma) = \pm
 \upsilon^2 \pi(\gamma')$, for some $\upsilon \in U$, then,
\[
\gamma / \overline{\gamma} = \pm \upsilon^2 \gamma'
/\overline{\gamma}' \qquad \mbox{and consequently,} \qquad
\upsilon^{-1} \gamma \overline{\gamma}' = \pm \upsilon \overline
\gamma \gamma'. 
\]
Since $N(\upsilon) = \pm 1$, and
 $\upsilon^{-1} = N(\upsilon) \overline{\upsilon}$, we have that either 
\[
\overline{\upsilon \overline \gamma \gamma'} = \upsilon
  \overline \gamma \gamma' \qquad \mbox{or} \qquad \overline{\upsilon
    \overline \gamma \gamma'} = -\upsilon \overline \gamma \gamma'.
\]
In the first case $\upsilon \overline{\gamma} \gamma' = n$ for some
rational integer $n$, and thus, $\upsilon N(\gamma) \gamma' = n
\gamma$, in which case $\gamma \sim \gamma'$.  In the second case,
$\upsilon \overline{\gamma} \gamma' = n \sqrt{d}$ for some rational
integer $n$, in which case $\upsilon N(\gamma) \sqrt{d} \gamma' = n d
\gamma$, which proves that $\gamma \sim \sqrt{d} \gamma'$.  
\end{proof}

\begin{lemma}
\label{lemma:3}
\begin{enumerate}
\item Given $\gamma \in \mathcal O$ there exists $\alpha \in \mathcal V$
such that $\gamma \sim \alpha$. 
\item If $\alpha, \alpha' \in \mathcal V$ are such that $\alpha \sim
  \alpha'$, then $| N(\alpha) | = | N(\alpha') |$. 
\end{enumerate}  
\end{lemma}
\begin{proof}
Suppose $\gamma = c + e \omega$ for $c, e \in \Z$.  Let $a =
c/\mathrm{gcd}(c,e)$ and $b = e/\mathrm{gcd}(c,e)$ and set $\alpha =
a + b \omega$.  Then clearly, $\alpha \in \mathcal V$ and since
$\gamma = \mathrm{gcd}(c,e) \alpha$, $\gamma \sim \alpha$.  

To prove 2, we note that if $\alpha \sim \alpha'$ then there exist integers $n$ and $m$ and a unit $\upsilon$ with
$\upsilon m \alpha' = n \alpha$.  With $\alpha = a + b \omega$ ($a$
and $b$ relatively prime integers) as before, we have
\[
\upsilon \alpha' = \frac{n a}{m} + \frac{n b}{m} \omega.
\]
Clearly then, since $na/m$ and $n b/m$ are relatively prime integers,
and $a$ and $b$ are relatively prime integers, we have that $n/m = \pm
1$. Consequently, there is an $\upsilon' \in U$ so that $\alpha = \upsilon' \alpha'$ and hence
$N( \alpha ) = \pm N(\alpha')$.  
\end{proof}

Our equivalence relation induces an equivalence relation on visible
points, and we will represent the corresponding equivalence classes by
$\widetilde{\mathcal V}$.  We will represent the equivalence class of
$\alpha \in \mathcal V$ by $\widetilde{\alpha}$.  The content of
Lemma~\ref{lemma:3} is that $\varphi \circ \pi$ is a two-to-one map of
$\widetilde{\mathcal V}$ onto $\mathbb S$ and the absolute value of
the norm is a well-defined function on $\widetilde{\mathcal V}$.  

Now, suppose $\xi \in \mathbb S$ with $\xi= \varphi(\beta)$ for some
$\beta \in \mathcal N$.  Suppose 
$\alpha, \alpha' \in \mathcal V$ are two non-equivalent visible points with 
\[
\varphi \circ \pi(\alpha) = \varphi \circ \pi(\alpha') = \xi.
\]
We define
\begin{equation}
\label{eq:14}
M(\xi) = \sqrt{|N(\alpha) N(\alpha')|},
\end{equation}
and use this to impose a partial ordering on $\varphi(\mathcal N)$ by
setting
\[
\mathcal M(t) = \{ \xi \in \varphi(\mathcal N) : M(\xi) \leq t \}
\]
and
\begin{equation}
\label{eq:10}
\mathcal M = \bigcup_{n=1}^{\infty} \mathcal M(n).
\end{equation}
\begin{thm}
\label{thm:2}
Let $K$ be a real quadratic extension of $\Q$.  Then $\mathcal M$
is equidistributed in $\mathbb S = \R / \Lambda(U^2).$ 
\end{thm}

This necessity to pass to equivalence classes is a consequence of the
infinitude of the unit group and suggests natural questions about
the distribution of the unit group itself.  The most natural question
is whether the unit group is equidistributed in $\R / \Z$.  
\begin{thm}\label{thm:4} Let $K$ be a real quadratic extension of $\Q$.  There is
  no ordering of the unit group so that it is equidistributed modulo
  1. That is, if $(\upsilon_n)_{n=1}^{\infty}$ is an ordering of $U$
  then $( \upsilon_n + \Z )_{n=1}^{\infty}$ is not equidistributed in
  $\R / \Z$.  
\end{thm}

\begin{proof}
It is enough to show that  $\left( \{ \upsilon_n \}
\right)_{n=1}^{\infty}$, the sequence of fractional parts of
$(\upsilon_n)_{n=1}^{\infty}$ converges  to 0. 
We may take $\epsilon$, a fundamental unit, to lie in $(0,1)$. The
unit group is $\{\pm \epsilon^n, \pm \overline{\epsilon}^n :n\in
\Z\}$.  Since $\epsilon \in (0,1)$, the powers of $\epsilon$ decrease
monotonically to 0. Therefore the sequence $-\epsilon^n$ modulo 1
increases monotonically to 0.  

First, we consider $d\not \equiv 1 \bmod 4$ so that 
$\epsilon^n=a_n+b_n\sqrt{d}$ where $a_n, b_n \in \Z$ and $ \overline{\epsilon}^n=a_n-b_n\sqrt{d}$.  Since $a_n\in \Z$, $\overline{\epsilon}^n\equiv - \epsilon^n \bmod 1$ and the powers of  $\overline{\epsilon}$ increases monotonically to 0 modulo 1.  Similarly, $-\overline{\epsilon}^n \equiv \epsilon^n \bmod 1$ so that these units decrease to 0. 

If $d\equiv 1 \bmod 4$, $\epsilon^n= a_n+b_n(1+\sqrt{d})/2$ for integers $a_n$ and $b_n$.  If $b_n$ is even, the above argument shows that the units are not equidistributed, so we will assume that $b_n$ is odd.  The conjugate $\overline{\epsilon}^n =a_n+b_n(1-\sqrt{d})/2$, and 
\[ \epsilon^n \equiv \frac12+\frac{b_n}2 \sqrt{d} \bmod 1,  \qquad  \overline{\epsilon}^n \equiv \frac12 - \frac{b_n}2 \sqrt{d} \bmod 1.\] 
By assumption, $\epsilon^n $ is positive and converges to 0.  Therefore $b_n\sqrt{d}/2 \rightarrow 1/2$ from the right. It follows that $-b_n\sqrt{d}/2 \rightarrow -1/2$ from the left.  As a consequence, $\overline{\epsilon}^n$ converges to 0 from the left. The other cases follow similarly.
\end{proof}

\section{An Outline of the Proofs}
\label{sec:an-outline-proof}

For the moment, and for concreteness we restrict our attention to
imaginary quadratics and the partially ordered set $\mathcal M$.  The characters of $\T$
are given by $\xi \mapsto \xi^m$ for $m \in
\Z$.  It is immediate that  \[
\int_{\T} \xi^m d\mu(\xi) = \left\{
\begin{array}{ll}
1 & \mbox{if } m = 0; \\
0 & \mbox{otherwise,}
\end{array}
\right.
\]
and thus, in order to use Weyl's criterion we need to show that 
\begin{equation}
\label{eq:2}
\lim_{n \rightarrow \infty} \frac{1}{\# \mathcal M(n)} \sum_{\beta
  \in \mathcal M(n)} \beta^m = 0  \qquad \mbox{for all} \qquad
m \neq 0.
\end{equation}
We set 
\[
a^{(m)}_n = \sum_{\beta \in \mathcal M(n) \setminus \mathcal M(n-1)} \beta^m 
\qquad
\mbox{and}
\qquad
A^{(m)}(x) = \sum_{n \leq x} a^{(m)}_{n} \qquad x > 0.
\]
Equation~(\ref{eq:2}) reduces to
\begin{equation}
\label{eq:4}
\lim_{x\rightarrow\infty} \frac{A^{(m)}(x)}{A^{(0)}(x)} = 0 \qquad
\mbox{for all} \qquad m \neq 0.
\end{equation}

\subsection{Serre's Idea}
\label{sec:serres-idea}

Following ideas of Serre in \cite{MR0263823}, we view 
$A^{(m)}(x)$ as the summatory function for the coefficients of an
$L$-series.  Analytic properties of the associated $L$-function can
then be used to determine the asymptotics of these summatory functions
as $x \rightarrow \infty$.   

Specifically, we define the $L$-series 
\[
\Xi^{(m)}(s) = \sum_{n=1}^{\infty} \frac{a^{(m)}_n}{n^{s}} =
\sum_{\beta \in \mathcal M} \frac{\beta^m}{N(\pi^{-1}(\beta))^s}
= \frac12 \sum_{\alpha \in \mathcal V} \frac{(\alpha/\overline{\alpha})^m }{N(\alpha)^s},
\]
and in particular, 
\[
\Xi^{(0)}(s) = \frac12 \sum_{\alpha \in \mathcal
  V} \frac{1}{N(\alpha)^s}. 
\]
The factor of $1/2$ is a result of the fact that restricted to $\mathcal V$, $\pi$ is two-to-one with 
$\pi^{-1}(\beta)= \pm \alpha$. 
Since $|N(n \alpha)| = n^2 |N(\alpha)|$ we find that
\begin{equation}
\label{eq:3}
\Xi^{(0)}(s) = \frac{1}{ 2 \zeta(2s)} \sum_{\alpha\in\mathcal V}
\sum_{n=1}^{\infty} \frac{1}{n^{2s} N(\alpha)^s} = \frac{1}{2 \zeta(2s)} \sum_{\alpha\in\mathcal V}
\sum_{n=1}^{\infty} \frac{1}{N(n \alpha)^s} 
\end{equation}
where $\zeta(s)$ is the Riemann zeta function.  Every non-zero in
$\mathcal O$ can be uniquely expressed as $n \alpha$ where $n$ is a
positive (rational) integer and $\alpha$ a visible point, so the
double sum in (\ref{eq:3}) can be replaced by a single sum over
$\mathcal O$.

If $\mf c = \gamma \mathcal O$ is a principal ideal, then
the norm of $\mf c$, $N \mf c$, is equal to $|N(\gamma)|$.  (The
absolute values are not necessary here, but they will be in the sequel
when we consider the real quadratic case, and at any rate, they do no
harm).  Two integers generate the same ideal if their quotient is a
unit, and thus, we may write $\Xi^{(0)}(s)$ as a sum over principal
ideals, so long as we compensate for the number of roots of unity in
$K$.  That is,

\[
\Xi^{(0)}(s) = \frac{w}{2\zeta(2s)} \sum_{\mf c \in I_0}
\frac{1}{(N \mf c)^s},
\]
where $I_0$ is the group of principal ideals in $\mathcal O$ and $w$
is the number of roots of unity in $K$ (4,6 or 2 depending 
on if $d = -1, -3$ or any other square-free negative number).

For a character $\chi$ defined on the group of principal ideals,
we define the {\em partial} $L$-series 
\[
L_{K,0}(s,\chi)= \sum_{\mf c \in
  I_0} \frac{\chi(\mf c)}{(N \mf c)^s}.
\]
If $\chi$ is the trivial character, then the resulting $L$-series,
\[
\zeta_{K,0}(s) = \sum_{\mf c \in I_0} \frac{1}{(N \mf c)^s}
\]
is the {\em partial} Dedekind zeta function of $K$ over the class of
principal ideals.  Therefore we can write 

\[
\Xi^{(0)}(s) =  \frac{w\zeta_{K,0}(s)}{2\zeta(2s)}.
\]

Let $\chi_k$ denote the principal Dirichlet character of modulus $k$,
defined by $\chi_k(\gamma)=1$ if $\gcd(N(\gamma),k)=1$ and 0 otherwise.
This character has a natural extension to a character on $I_0$, which
we will denote by $\chi_k$ as well.  Specifically, $\chi_k(\mf c) =
\chi_k(\gamma)$ if $\mf c=\gamma \mathcal O$.  

We will also write
\[ \zeta(s, \chi_k) = \sum_{n=1}^{\infty} \frac{\chi_k(n)}{n^s} \]
where, in this context, we consider $\chi_k$ to be the principal
Dirichlet character of modulus $k$ on $\Z$.  We will require these more general
$L$-series in the proof of Theorem \ref{thm:1} and \ref{thm:2}.

The following lemma is more robust than we need at present  to relate
$\Xi^{(m)}(s)$ to partial $L$-series, but we will require this
stronger form in the sequel.  Let $\xi$ be a generator of the unit
group and define 
\[ 
S(\chi)= \sum_{\ell=1}^w \left[ \chi(\xi)\right]^{\ell}.
\]
We will call a character $\chi$ on $\mathcal O$ {\em admissible} if
$S(\chi)=0$ or $\chi$ is well-defined on $\widetilde{\mathcal O}$. If
$\chi$ is admissible, define $\chi^*$ to be identically zero if $S=0$
and otherwise define the character $\chi^*:I_0\rightarrow \T$ by
 \[ \chi^*(\mf{c})=\chi(\gamma) \qquad \mf{c}=\gamma \mathcal{O}.\]

\begin{lemma}
\label{lemma:2}   If $K$ is an imaginary quadratic extension of $\Q$ and $\chi$ is an admissible character, 
then 
\[
 \sum_{\alpha \in \mathcal V} \frac{ \chi(\alpha) \chi_k(\alpha)}{N(\alpha)^s}=S(\chi)\:
\frac{ L_{K,0}(s, \chi^* \chi_k)}{\zeta(2s,\chi_k)}.
\]


\end{lemma}

\begin{proof}

We can perform  maneuvers similar to the $\Xi^{(0)}$ case to write 
\[
L(s) = \frac{1}{\zeta(2s,\chi_k)} \sum_{\alpha \in \mathcal V}
\sum_{n=1}^{\infty} \frac{\chi(\alpha)\chi_k(\alpha) \chi_k(n)}{n^{2s} N(\alpha)^s}
\] where $L(s)$ is the left hand $L$-series. 
The unit group $U$ consists of the roots of unity in $K$ with generator $\xi$, and 
$U$ acts on $\mathcal V$ by the map $\alpha \mapsto \xi \alpha$.  Let
$\mathcal V'$ be a complete collection of representatives of the
orbits of this action.  Then, 
\begin{align*}
L(s) &= \frac{1}{\zeta(2 s,\chi_k)} \sum_{\alpha \in \mathcal V'}
\sum_{n=1}^{\infty} \sum_{\ell=1}^w
\frac{\chi(\xi)^{\ell} \chi(\alpha) \chi_k(\alpha) \chi_k(n)  }{n^{2s}
  N(\alpha)^s}.\end{align*}
  Since $\chi_k(n)\chi_k(\alpha)=\chi_k(n\alpha)$, 
 \[  L(s) = \frac{1}{\zeta(2 s,\chi_k)}  \bigg[ \sum_{\ell=1}^w
\chi(\xi)^{\ell} \bigg] \sum_{\alpha \in \mathcal V'}
\sum_{n=1}^{\infty}
\frac{\chi(\alpha)\chi_k(n\alpha)}{n^{2s} N(\alpha)^s}. \]
We may now assume that  $\chi$ is well-defined on $\widetilde{\mathcal
  O}$, so that  
 $\chi(\alpha) = \chi(n\alpha)$.  It follows that 
\[
L(s) = \frac{S}{\zeta(2 s,\chi_k)} \sum_{\gamma \in \mathcal O}
\frac{\chi(\gamma) \chi_k(\gamma)}{N(\gamma)^s} = \frac{S}{\zeta(2s,\chi_k)}
\sum_{\mf c \in I_0} 
\frac{\chi^*(\mf c)\chi_k(\mf c) }{(N \mf c)^s}. \qedhere
\]
\end{proof}

Returning to $\Xi^{(m)}$,  define the character $\tau^{(m)}: I_0 \rightarrow \T$ on
the group of principal ideals by setting 
\[
\tau^{(m)}(\mf c) = (\gamma/\overline{\gamma})^m \qquad \mf c = \gamma \mathcal O.
\]
This character is admissible.  Moreover,
\[
 S=\sum_{\ell=1}^w (\xi/{\overline \xi})^{\ell m} =  \sum_{\ell=1}^w
\left( \xi^{2m} \right)^{\ell} = \left\{
\begin{array}{ll}
w & \mbox{if } 2 m \equiv 0 \bmod w; \\ & \\
0 & \mbox{otherwise.}
\end{array}
\right.
\]
We apply the lemma  using $\chi_r=\chi_0$ and $\chi^*=\tau^{(m)}$. We conclude that
if $2m \not \equiv 0
\bmod w$  the summatory function, $A^{(m)}(x)$, is identically 0,
and therefore (\ref{eq:4}) is trivially satisfied.  That is, it
suffices to establish (\ref{eq:4}) for $2 m \equiv 0 \bmod w$.  
 It follows that we may write
\[
\qquad \Xi^{(m)}(s) = \frac{wL_{K,0}(s, \tau^{(m)})}{2\zeta(2s)}.
\]

\subsection{The Wiener-Ikehara Tauberian Theorem}

The sum defining $\zeta_{K,0}(s)$ is absolutely convergent in the
half-plane $\rp{s} > 1$.  It follows, therefore that $L_{K,0}(s,\tau^{(m)})$ is
absolutely convergent in the same half-plane, and by Morera's Theorem, 
the resulting function of $s$ is analytic.  Via the following
corollary to the Wiener-Ikehara Tauberian Theorem, an analytic (or
meromorphic) continuation of  $\Xi^{(m)}(s)$ to a half-plane containing
the line $s=1$ yields asymptotic information about the summatory
function $A^{(m)}(x)$. 
\begin{thm}[Ikehara, Wiener]
Suppose 
\[
f(s) = \sum_{n=1}^{\infty} \frac{a_n}{n^s}; \qquad \qquad a_n \geq 0
\]
converges to an analytic function of $s$ on the half-plane $\rp{s} >
1$, and
\[
g(s) = \sum_{n=1}^{\infty} \frac{b_n}{n^s}; \qquad  \qquad b_n \in \C
\] 
is such that $| b_n | \leq a_n$.  If $f(s)$ has a meromorphic continuation to
a neighborhood of the line $\rp{s} \geq 1$ with a single simple pole at $s=1$
and $g(s)$ has an analytic continuation to a neighborhood of the line
$\rp{s} \geq 1$, then 
\[
\sum_{n \leq x} a_n \sim r x; \qquad \mbox{where} \qquad r =
\mathop{\Res}_{s=1} f(s), 
\]
and
\[
\sum_{n \leq x} b_n = o(x).
\]
\end{thm}
The proof of this formulation of the Wiener-Ikehara Tauberian Theorem
can be found in \cite[Ch. VIII, \S 3]{MR0160763}.

\subsection{Hecke's Functional Equation}

It remains to show that $\Xi^{(0)}(s)$ has a meromorphic continuation
to a neighborhood of $\rp{s} \geq 1$ with a single simple pole at
$s=1$ and for all $m\neq 0$, $\Xi^{(m)}(s)$ has an analytic
continuation to a neighborhood of $\rp{s} \geq 1$.  Then, by the
Tauberian theorem, (\ref{eq:4}) will have been met
with  $g(s)=\Xi^{(m)}(s)$  and $f(s)=\Xi^{(0)}(s)$, and we will have
arrived at our equidistribution theorem.  Since $1/ \zeta(2s)$ is
analytic in the half-plane $\rp{s} > 1/2$, it suffices to find a
meromorphic (respectively analytic) continuation for $\zeta_{K,0}(s)$
and $L_{K,0}(s,\tau^{(m)})$.

In \cite{MR1544392}, Hecke produced a functional equation for a large
class of $L$-functions, 
which includes $\zeta_{K,0}(s)$ and $L_{K,0}(s,\tau^{(m)})$.  Since we
will ultimately need the analytic continuation for more complicated
$L$-functions than $L_{K,0}(s,\tau^{(m)})$, we review some of Hecke's
theory.  

Here we may suppose that $K$ is an arbitrary number field.  We denote
the group of fractional ideals of $K$ by $\wt I$ and the subgroup of
principal fractional ideals by $\wt I_0$.  Given a set of fractional
ideals $\wt J$, we will represent its subset of integral ideals by
$J$.  Given a (proper) integral ideal $\mf m \in I$, we define the
subgroups of $\wt I$ given by
\[
\wt I_{\mf m} = \{ \mf a \in I : \mf a = \mf b \mf c^{-1} \mbox{ with }
(\mf{b c}, \mf m) = \mathcal O \} \qquad \mbox{and} \qquad \wt I_{\mf m,0}
= \wt I_{\mf m} \cap \wt I_0. 
\]
That is $\wt I_{\mf m,0}$ is exactly the subset of principal fractional
ideals in $\wt I_{\mf m}$.

Notice that
\[
(\mathcal O/\mf m)^{\times} = \{ \gamma + \mf m : \gamma \mathcal O
\in I_{\mf m,0} \} 
\]
is a group under multiplication.  A character $\chi: (\mathcal O/\mf
m)^{\times} \rightarrow \T$ is called a Dirichlet character of modulus 
$\mf m$.  As with the usual Dirichlet characters defined on the
rational integers, we extend the definition of $\chi$ to all of
$\mathcal O$ by setting
\[
\chi(\gamma) = \left\{
\begin{array}{ll}
\chi(\gamma + \mf m) & \quad \mbox{if } \gamma \mathcal O \in I_{\mf m,0}; \\ 
0                    & \quad \mbox{otherwise}. 
\end{array}
\right.
\]
The character which is identically one on $I_{m,0}$ is called the
principal Dirichlet character of modulus $\mf m$ and denoted
$\chi_{\mf m}$.  The Dirichlet character $\chi_k$ introduced in
Section~\ref{sec:serres-idea} is equal to the Dirichlet character
$\chi_{k \Ox}$. 

If $K$ has $r_1$ real places and $r_2$ complex conjugate places,
 then we may construct a character $\chi_{\infty}$ on
$K^{\times}$ by pulling back a character $(\R^{\times})^{r_1} \times
(\C^{\times})^{r_2} \rightarrow \T$ through the natural embedding $K
\hookrightarrow \R^{r_1} \times \C^{r_2}$.  Such a character
is called an {\em Archimedean} character of $K$.

A Hecke character $\psi$ is a character on $\wt I$ whose restriction to
principal integral ideals, can be written as 
\begin{equation}
\label{eq:5}
\psi\big(\mf c \big) = \chi(\gamma) \chi_{\infty}(\gamma) \qquad \mf c
= \gamma \Ox,
\end{equation}
where $\chi$ is a Dirichlet character for some ideal $\mf m$ and
$\chi_{\infty}$ is an Archimedean character of $K$.
Equation~(\ref{eq:5}) implies the {\em unit consistency condition}:
If $\upsilon \in U$ then $\chi(\upsilon) \chi_{\infty}(\upsilon) =
1$. In fact, it can be shown that any pair of characters consisting of
a Dirichlet character and an Archimedean character which satisfies
the unit consistency condition can be derived 
from a Hecke character.  A {\em principal} Hecke character is one for which
$\chi$ is a principal Dirichlet character for some modulus, and
$\chi_{\infty}$ is identically one.  The trivial Hecke character,
which is identically one on the set of non-zero ideals in $\wt I$, is
an example of a principal character.

To each Hecke character we associate the $L$-series
\[
L_K(s,\psi) = \sum_{\mf c \in I}
\frac{\psi(\mf c)}{N(\mf c)^s} \qquad \mbox{and} \qquad
L_{K,0}(s,\psi) = \sum_{\mf c \in I_0} \frac{\psi(\mf c)}{N(\mf
  c)^s}.
\]
The former is an example of a Hecke $L$-series, the latter a partial
Hecke $L$-series over the principal ideal class.  The sums defining
these functions are absolutely convergent in the half-plane $\rp{s} >
1$ and therefore in that region they are analytic functions of $s$.
The Hecke $L$-series corresponding to the trivial character define
respectively the Dedekind zeta function of $K$ and the partial
Dedekind zeta function over the principal ideal class of $K$.  

Hecke demonstrated that his $L$-functions and partial $L$-functions
have functional equations analogous to that for the Riemann zeta
function.  However, since we will only be interested in the analytic
properties of $L_K(s, \psi)$ and $L_{K,0}(s, \psi)$ in a neighborhood of the
half plane $\rp{s} \geq 1$ (and even then for only a very limited set of
Hecke $L$-series), and since it simplifies the exposition, we will be
satisfied with the following corollary to Hecke's Theorem.
\begin{thm}[Hecke]\label{thm:hecke}
If $\chi$ is a principal Hecke character, then $L_K(s, \chi)$ and
$L_{K,0}(s,\chi)$ have meromorphic continuations to all of $\C$ with a
single simple pole at $s=1$.  If $\chi$ is a non-principal character,
then $L_K(s,\chi)$ and $L_{K,0}(s,\chi)$ have analytic continuations to
entire functions of $s$.   
\end{thm}

\subsection{The Proof of Theorem~\ref{thm:3} }\label{section:thm3proof}

We finally have all the tools necessary to prove Theorem~\ref{thm:3}.
Theorem~\ref{thm:hecke} implies that $\zeta_{K,0}(s)$ has a meromorphic continuation to $\C$ with the exception of a single simple pole at $s=1$. Therefore $\Xi^{(0)}(s) $ has a meromorphic continuation to a neighborhood of $\rp{s}\geq 1$ with the exception of a single simple pole at $s=1$.
As a result, it suffices to show that, if $2m \equiv 0 \bmod w$ is a non-zero
integer, then $\tau^{(m)}$ is a non-principal Hecke character.  Using the Tauberian theorem this will imply
the analytic properties of $\Xi^{(m)}(s)$ necessary to show that the
summatory function $A^{(m)}$ satisfies (\ref{eq:4}) which in turn 
implies that $\mathcal M$ satisfies Weyl's Criterion and hence is
equidistributed.  We will prove a slightly more general result, which we will use in the proof of Theorem \ref{thm:1}.

\begin{lemma}\label{lemma:5} If $2m \equiv 0 \bmod w$ then 
the character $\mf c \mapsto \chi_k(\mf c) \tau^{(m)}(\mf c)$ is a
Hecke character and is principal only for $m=0$.  
\end{lemma}

\begin{proof}
The map from $\C^* \rightarrow \T $ defined by $z\to
(z/\overline{z})^m$ is an Archimedian character, and induces a map
$\chi_{\infty}$ on $\mathcal O$. For all $\mf c = \gamma
\mathcal O \in I_0$,  
\[ 
\chi_k( \mf{c}) \tau^{(m)}(\mf c) = \chi_k(\gamma)\chi_{\infty}(\gamma).
\] 
If $2m \equiv 0 \bmod w$, and $\upsilon$ is any unit, then 
\[
\chi_k(\upsilon) \chi_{\infty}(\upsilon) = \chi_{\infty}(\upsilon) =
\frac{\upsilon^{2m}}{N(\upsilon)^m} = 1,
\]
since $\upsilon$ must be a $w$th root of unity.  That is, $(\chi_k,
\chi_{\infty})$ satisfy the unit consistency condition, and $\chi_k\tau^{(m)}$
is the restriction of some Hecke character on $I_0$.  This character
is clearly non-principal if $m \neq 0$.
\end{proof}

\section{The Proof of Theorem~\ref{thm:2}}\label{section:realproof}

We now turn to the equidistribution of $\mathcal M$ for real quadratic
extensions of $\Q$.  The characters of $\mathbb{S}$ are of the form
\[
\tau^{(m)}: \xi \mapsto \exp\{2 \pi i m \xi/ \log| \epsilon^2 |\}.
\]
We shall reuse the notation $\tau^{(m)}$, $a_n^{(m)}$, $\Xi^{(m)}(s)$
etc.~for objects in the real quadratic case which are analogous to the
complex quadratic case.  

In order to use Weyl's criterion we must show that
\[
\lim_{n \rightarrow \infty} \frac{1}{\# \mathcal M(n)} \sum_{\xi \in
  \mathcal M(n)} \exp\{2 \pi i m \xi/ \log| \epsilon^2 |\} = 0 \qquad
\mbox{for all} \qquad m \neq 0. 
\]
As in the imaginary quadratic case, this can be deduced from the
analytic properties of 
\[
\Xi^{(m)}(s) = \sum_{n=1}^{\infty} \frac{a_n^{(m)}}{n^s},
\]
where in the case at hand,
\[
a_n^{(m)} = \sum_{\xi \in \mathcal M(n) \setminus \mathcal M(n-1)}
\tau^{(m)}(\xi). 
\]
If $\xi \in \mathcal M(n) \setminus \mathcal M(n-1)$, then $n-1 <
M(\xi) \leq n$.  As we shall see, $M(\xi)$ is always an integer, and
thus we may write
\begin{equation}
\label{eq:13}
\Xi^{(m)}(s) = \sum_{\xi \in \mathcal M}
\frac{\tau^{(m)}(\xi)}{M(\xi)^s}.
\end{equation}
The equidistribution theorem for real quadratic fields will be proved
if we can show that  $\Xi^{(m)}(s)$ is analytic in a
neighborhood of $\rp{s} \geq 1$, except in the case $m=0$ where it has
a single simple pole at $s=1$.

\begin{lemma}
 If $\xi \in \mathcal M$ is such
that 
\[
\xi = \varphi(\pi(\alpha)) \qquad \mbox{for} \qquad \alpha \in
\mathcal V,
\]
then, 
\[
M(\xi) =  \frac{\sqrt{d}\:  |N(\alpha)|}{\mathrm{gcd}(N(\alpha), d)}.
\]
\end{lemma}
\begin{proof} First, we will prove the lemma in the case where $d\not \equiv 1 \bmod 4$.  
There exist relatively prime rational integers $a$ and $b$ so that
$\alpha = a + b \sqrt{d}$.  We may construct the non-equivalent
visible point $\alpha'$ by setting 
\[
\alpha' = \frac{b d}{\mathrm{gcd}(a, d)} + \frac{a}{\mathrm{gcd}(a,
  d)} \sqrt{d}.
\]
It follows that
\[
M(\xi) = \sqrt{| N(\alpha) N(\alpha') |} =
\frac{\sqrt{d} \: |N(\alpha)|}{\mathrm{gcd}(a, d)}.
\]
Since $N(\alpha) = a^2 - b^2 d$ it is easy to see that $\gcd(a, d) =
\gcd(N(\alpha), d)$.

If $d \equiv 1  \bmod 4$, then there are relatively prime integers $a$ and $b$ so that $\alpha = a+b\omega$.  We have
\[ \sqrt{d} \alpha = b ( d-1)/2 -a + (2a+b)\omega. \] It is elementary to see that 
\[ \gcd \big( b ( d-1)/2 -a, 2a+b \big) = \gcd(N(\alpha), d)\]
so that the non-equivalent visible point $\alpha'$ is 
\[ \alpha' = \frac{ b ( d-1)/2-a + (2a+b)\omega}{\gcd(N(\alpha), d)  }. \]
It follows that 
\[ M(\xi) = \sqrt{|N(\alpha)N(\alpha')|} = \frac{\sqrt{d} |N(\alpha)|}{\gcd(N(\alpha),d)}\]
in this case as well. 
\end{proof}

The two-to-one map from $\widetilde{\mathcal V}$ to $\mathcal M$
allows us to write (\ref{eq:13}) as
\[
\Xi^{(m)}(s) = \frac{ 1 }{2d^{\frac{s}2}}\sum_{\widetilde{\alpha} \in
  \widetilde{\mathcal V}} \gcd(N(\alpha), d)^s \frac{\tau^{(m)}\big(
  \varphi \circ \pi(\alpha) \big) }{| N(\alpha) |^s}
\]
since the summand is independent of the representative $\alpha \in \widetilde{\alpha}$.
Given $r$ dividing $d$,
let 
\[
\widetilde{\mathcal V}(r) = \{ \widetilde \alpha \in \mathcal V :
\gcd(N(\alpha), d) = r \}.
\]
Using this notation,
\[
\Xi^{(m)}(s) =\frac{ 1 }{2d^{\frac{s}2}} \sum_{r | d} r^s \sum_{\widetilde{\alpha} \in
  \widetilde{\mathcal V}(r)}  \frac{\tau^{(m)}\big(
  \varphi \circ \pi(\alpha) \big) }{| N(\alpha) |^s}. \]

Define the character $\eta_r$ on $\mathcal O$ as
\[ \eta_r(\gamma)= \left\{ \begin{aligned} 
1 & \quad \mbox{if} \ (N(\gamma), d) = r\\ 0 & \quad \mbox{otherwise.}\end{aligned} \right. \]
This is well-defined on $\widetilde{\alpha}$, so
\[
\Xi^{(m)}(s) = \frac{ 1 }{2d^{\frac{s}2}}\sum_{r | d} r^s \sum_{\widetilde{\alpha} \in
  \widetilde{\mathcal V}}  \frac{\eta_r(\alpha) \tau^{(m)}\big(
  \varphi \circ \pi(\alpha) \big) }{| N(\alpha) |^s}. \]

The character $\eta_r$ can be extended to the non-zero principal ideals, $I_0$, by defining $\eta_r(\mf{c})=\eta_r(\gamma)$ if $\mf c = \gamma \mathcal O$.  
Let $\chi_k$ be the principal character of modulus $k$. 
We can rewrite $\eta_r$ in terms of these principal characters.
\begin{lemma}\label{lemma:7}  The character $\eta_r$ satisfies 
\[ \eta_r = \sum_{\ell | d/r} \mu(\ell) \sum_{k|r\ell} \mu(k) \chi_k.\]
\end{lemma}

\begin{proof} 
For square-free positive integers $s$, the character $\omega_s (\gamma)= \sum_{k|s} \mu(k) \chi_k (\gamma)$ is 1 if $(N(\gamma),s)=s$ and 0 otherwise. It follows that $\sum_{\ell | d/r} \mu(\ell) \omega_{r\ell}$
defines the character $\eta_r$. 
\end{proof}

Using this lemma
\[ 
\Xi^{(m)}(s)   = \frac{ 1 }{2d^{\frac{s}2}}
\sum_{r | d} \sum_{\ell | d/r}  \sum_{k|r\ell}  \big(r^s \mu(\ell) \mu(k)\big)
\sum_{\widetilde{\alpha} \in \widetilde{\mathcal V} } \frac{\chi_k(\alpha) \tau^{(m)}\big(
  \varphi \circ \pi(\alpha) \big) }{| N( \alpha) |^s}. 
\]

We now prove two lemmas that will enable us to complete the proof of Theorem~\ref{thm:2}.
\begin{lemma}\label{lemma:8} If $K$ is a real quadratic number field and $\chi$ is an admissible character, then 
\[ \sum_{\widetilde{\alpha} \in \widetilde{\mathcal V}} \frac{ \chi({\alpha})  \chi_k({\alpha})}{|N({\alpha})|^s}=  \frac{ L_{K,0}(s, \chi^* \chi_k)}{\zeta(2s,\chi_k)}.\]
\end{lemma}
Note that if $\chi$ is an admissible character, then the terms are independent of  the representative of $\widetilde{\alpha} \in \widetilde{\mathcal V}$.

\begin{proof}
The character  $\chi_k(\widetilde{\alpha})$ is well defined in $\widetilde{\mathcal V}$ by Lemma  \ref{lemma:3}.
and since $\chi$ is well defined on equivalence classes,
\begin{align*}
 \sum_{\widetilde{\alpha} \in \widetilde{\mathcal V}} \frac{ \chi({\alpha})  \chi_k({\alpha})}{|N({\alpha})|^s} 
 &= \frac{1}{\zeta(2 s,\chi_k)} \sum_{\widetilde{\alpha} \in {\widetilde{\mathcal V}}}
\sum_{n=1}^{\infty}
\frac{ \chi({\alpha}) \chi_k({\alpha}) \chi_k(n)  }{n^{2s}
  |N({\alpha})|^s}  \\ 
 &= \frac{1}{\zeta(2 s,\chi_k)} \sum_{\widetilde{\gamma} \in {\widetilde{\mathcal O}}}
\frac{ \chi({\gamma})\chi_k({\gamma})}{|N({\gamma})|^s}\\
& = \frac{1}{\zeta(2 s,\chi_k)} \sum_{\mf c \in I_0} \frac{\chi^*(\mf c) \chi_k(\mf c)}{(N \mf c)^s}.
\end{align*}

\end{proof}

\begin{lemma}\label{lemma:6}
The character $\chi_k\: \tau^{(m)}(\varphi \circ \pi)$ is well-defined on $\widetilde{\mathcal O}$ and extends to a Hecke character which is principal only for $m=0$. 
\end{lemma}

\begin{proof} Let $\widetilde{\alpha} \in \widetilde{\mathcal V}$ and $\alpha \in \widetilde{\alpha}$.
The map $\pi$ satisfies $\pi(\alpha)=\pi(n\alpha)$ for all  $n\in \Z$, and $\varphi \circ
\pi(u\alpha)=\varphi \circ \pi (\alpha)$.  Therefore, since $\chi_k$
is well-defined on $\widetilde{\mathcal O}$, $\chi_r \:
\tau^{(m)}(\varphi \circ \pi)$ is as well.  The map
$\chi_{\infty}:\R^{\times} \times \R^{\times} \rightarrow \T$ defined
by
\[ 
(x,y) \mapsto \exp\{ 2 \pi i m \log|x/y|/\log |\epsilon|^2\} 
\] is
an Archimedean character which extends to $\tau^{(m)}$ in $K$.  As
$\chi_k$ is a Dirichlet character we need only verify the unit
consistency condition.  Let $\upsilon \in U$ be a unit, then
\[ 
\chi_k(\upsilon)\tau^{(m)}(\upsilon) = \chi_k(\upsilon)
\chi_{\infty}(\upsilon) = 1
\] 
since $\log (\upsilon/\bar{\upsilon}) = 0$.  The character is clearly
only principal when $m=0$.
\end{proof}

Lemma~\ref{lemma:6} proves that $\chi_k \tau^{(m)}$ is an admissible
character. Using Lemma~\ref{lemma:8} we conclude that
\[\sum_{\widetilde{\alpha} \in \widetilde{\mathcal V} }
\frac{\chi_k(\alpha) \tau^{(m)}\big( \varphi \circ \pi(\alpha) \big)
}{| N( \alpha) |^s} = \frac{L_{K,0}(s,\tau^{(m)}
\chi_k)}{\zeta(2s,\chi_k)}.\] Therefore,
\[ \Xi^{(m)}(s) = \frac{ 1 }{2d^{\frac{s}2}} \sum_{r | d} \sum_{\ell |
d/r} \sum_{k|r\ell} \big(r^s \mu(\ell) \mu(k)\big)
\frac{L_{K,0}(s,\tau^{(m)} \chi_k)}{\zeta(2s,\chi_k)}.\] That is,
$\Xi^{(m)}(s)$ is a finite sum of partial $L$-series of the form
$L_{K,0}(s, \tau^{(m)} \chi_k)$ divided by a finite number of factors
of the form $\zeta(2s,\chi_k)$. By Lemma \ref{lemma:6}, $\tau^{(m)}
\chi_k$ is a Hecke character and is principal only for
$m=0$. Therefore $L_{K,0}(s, \tau^{(m)}\chi_k)$ has an analytic
continuation to a neighborhood of the half-plane $\rp s \geq 1$ if $m\neq 0$ and a
meromorphic continuation to a neighborhood of the half-plane $\rp s \geq 1$ if $m=0$
with a simple pole at $s=1$.  The same is necessarily true of
$\Xi^{(m)}(s)$.  We can now use Weyl's criterion to prove Theorem
\ref{thm:3} invoking the Tauberian theorem to estimate the coefficient
sums.

\section{The Proof of Theorem \ref{thm:1}}

To prove the equidistribution of the set $\mathcal H$ we again use the
Tauberian theorem in conjunction with Weyl's criterion. We define the
$L$-series
\[ \Upsilon^{(m)} (s) = \sum_{\beta \in \mathcal H} \frac{\beta^m
}{H(\beta)^s}. \] As in the proof of Theorem \ref{thm:3} it is
sufficient to show that $\Upsilon^{(0)}(s)$ has a meromorphic
continuation to a neighborhood of $\rp{s}\geq 1 $ with a single simple
pole at $s=0$ and for all $m\neq 0$, $\Upsilon^{(m)}(s)$ has an
analytic continuation to a neighborhood of $\rp{s}\geq 1$.  We will do
so by rewriting $\Upsilon^{(m)}(s)$ in terms of other $L$-series whose
analytic properties are known.

\begin{lemma}\label{lemma:4} Suppose $\beta \in \mathcal H$ and
$\alpha \in \mathcal V$ with $\pi(\alpha)=\beta$. Then \[ H(\beta)^2 =
N(\alpha)/M(\alpha) \] where $M(\alpha)= \gcd(N(\alpha),2d)$.
\end{lemma}
Before proving this, we introduce the Mahler measure of integer
polynomials and its relationship with the absolute Weil height.  
Given a polynomial which factors over $\C$ as 
\[
f(z) = c \prod_{n=1}^N (z - \beta_n). 
\]
The {\em Mahler measure} of $f$ is defined to be 
\[
\mu(f) = |c| \prod_{n=1}^N \max\{1, |\beta_n| \}.
\]
If $f$ is irreducible and $f(\beta) = 0$ (that is, $\beta =
\beta_n$ for some $1 \leq n \leq N$), then the absolute Weil height
of $\beta$ is related to the Mahler measure of $f$ by the
following identity: 
\[
H(\beta)^N = \mu(f).
\]
(See, for instance \cite[Prop. 1.6.6]{MR2216774} for a proof.)
\begin{proof}[Proof of Lemma~\ref{lemma:4}]
Given $\beta = \alpha/\overline{\alpha} \neq \pm 1$ with $\alpha = a +
b \omega \in \mathcal V$, the minimal polynomial of $\beta$ over $\Q$
is given by 
\[
(z - \alpha/\overline{\alpha})(z - \overline{\alpha}/\alpha) = z^2 -
\left(\frac{\alpha^2 + \overline{\alpha}^2}{N(\alpha)} \right)z + 1 =
z^2 - \frac{2 a^2 + 2 a b (\omega + \overline{\omega}) + b^2(\omega^2
  + \overline{\omega}^2)}{a^2 + ab(\omega +
    \overline{\omega}) + b^2 \omega \overline{\omega}} z + 1.
\]
If we set 
\[
M = \gcd\big(a^2 + ab(\omega + \overline{\omega}) + b^2 \omega
\overline{\omega}, 2 a^2 + 2 a b (\omega + \overline{\omega}) +
b^2(\omega^2 + \overline{\omega}^2) \big),
\]
then the integer minimal polynomial of $\beta$ is given by
\[
f(z) = \frac{N(\alpha)}{M} z^2 - \frac{2 a^2 + 2 a b (\omega +
  \overline{\omega}) + b^2(\omega^2 
  + \overline{\omega}^2)}{M} z +  \frac{N(\alpha)}{M}. 
\]
Since $\beta$ and $\overline{\beta}$ are both on the unit circle, we
have 
\[
H(\beta)^2 = \mu(f) = \frac{N(\alpha)}{M}.
\]
It remains to show that $M = M(\alpha)$ as given in the statement of
the lemma. 

Case $d \equiv 1 \bmod 4$: In this situation, $\omega +
\overline{\omega} = 1$ and $\omega \overline{\omega} = (1 - d)/4$, and
thus,
\begin{equation}
\label{eq:6}
M = \gcd\left( a^2 + ab + b^2 \left(\frac{1 - d}{4}\right), 2 a^2 +
  2 ab + b^2 \left(\frac{1 + d}{2}\right) \right).
\end{equation}

We observe that right term minus twice the left term is $b^2d$. The
observation implies that if $p^n$ divides $M$, then $p$ divides $d$ or
$b$.

If $p$ is an odd prime and $p$ divides $b$ then since $p$
divides $N(\alpha)$, $p$ must also divide $a$, which cannot occur as
$a$ and $b$ are relatively prime; so $p$ divides $d$.  (Therefore the
odd part of $M$ must be square-free.)  If $p$ divides
$\gcd(N(\alpha),d)$ then the observation implies that $p$ divides
$M$ as well. 

If 2 divides $M$, then $b$ must be odd and $(1-d)/4$ must be even,
so that $d\equiv 1 \bmod 8$.  As a result, 4 cannot divide the
right hand term of (\ref{eq:6}).  Therefore, 2 divides both
$N(\alpha) $ and $2d$. Conversely, if 2 divides both $N(\alpha)$ and
$2d$ then 2 must divide the $M$.  Since $d$ is square-free,
$M=\gcd(N(\alpha), 2d)$.  

Case $d \not\equiv 1 \bmod 4$:  Here, $\omega + \overline{\omega} = 0$
and $\omega \overline{\omega} = -d$, and thus,
\[
M = \gcd ( a^2 - b^2 d, 2 a^2 + 2 b^2 d).
\]
We observe that twice the left term plus the right term is $4a^2$ and
the right term minus twice the left term is $4b^2d$.  Therefore, if
$p$ is an odd prime and $p^n$ divides the $M$, then $p^n$ divides
$a^2$ and $b^2d$. Since $a$ and $b$ are relatively prime, $p^n$
divides $a$ and $d$.  It follows that the odd part of the $M$ is
square-free.  This observation also demonstrates that if the odd prime
$p$ divides $d$ and $N(\alpha)$, then it divides $M$. 

Now we consider even divisors. First, assume that $d\equiv 2 \bmod
4$. If 2 divides $M$, then $a$ must be even and hence $b$ must be odd
as $\gcd(a,b)=1$.  Therefore $b^2d$ is congruent to 2 modulo 4, and it
follows that 4 cannot divide $M$. We conclude that $M$ is equal  to
$\gcd(N(\alpha), d)$.  Since $d\equiv 2 \bmod 4$, $N(\alpha) $ is
either odd or not divisible by 4, so that this   is equivalent to
$\gcd(N(\alpha), 2d)$.  

 Finally, assume that $d\equiv 3 \bmod 4$.  If 2 divides $M$ then 2
 divides $a^2-b^2d$ so that both $a$ and $b$ must be odd as they are
 relatively prime.  As a result,  the right term is not divisible by
 4. We conclude  that $M= \gcd(N(\alpha), 2d)$ in this case as
 well. 
\end{proof}

For each positive integer $r$ dividing $2d$, we define a character
$\eta_r$ on $\mathcal O$ as
\[ \eta_r(\gamma)= \left\{ \begin{aligned} 1 & \quad \mbox{if} \
(N(\gamma), 2d) = r\\ 0 & \quad \mbox{otherwise}.
\end{aligned} \right. \] This is identical to the definition of
$\eta_r$ in Section \ref{section:realproof}, but with $2d$ replacing
$d$ so that Lemma~\ref{lemma:7} applies. (In the $d\equiv 2 \bmod 4$ case, the gcd is square-free even though $2d$ is not, so the proof of Lemma~\ref{lemma:7} is valid in this case as well. )

We write \[ M(\alpha) = (N(\alpha),2d) = r \eta_r(\alpha)\] where
$r=(N(\alpha),2d).$ It follows that $M(\alpha) = \sum_{r|2d}
r\eta_r(\alpha)$ and since exactly one summand is $r$ and the rest
are zero,
\[ M(\alpha)^s = \sum_{r|2d}  r^s \eta_r(\alpha).\]

We can now rewrite $\Upsilon^{(m)}(2s)$ as
\[ \begin{aligned} \Upsilon^{(m)}(2s) = \frac12 \sum_{\alpha \in
\mathcal V} \frac{(\alpha/\overline{\alpha})^m
M(\alpha)^s}{N(\alpha)^s} = \frac12  \sum_{r|2d} r^s \sum_{\alpha
\in \mathcal V} \frac{(\alpha/\overline{\alpha})^m \eta_r(\alpha)
}{N(\alpha)^s}.
\end{aligned} \] 
It will suffice to show the required analytic properties of $\Upsilon^{(m)}(2s)$.
Using Lemma \ref{lemma:7} we can rewrite this as
\[\Upsilon^{(m)}(2s) = \frac12  \sum_{r|2d} \sum_{\ell | 2d/r}
\sum_{k|r\ell} \big(r^s \mu(\ell) \mu(k) \big) \sum_{\alpha \in
\mathcal V} \frac{(\alpha/\overline{\alpha})^m \chi_k(\alpha)
}{N(\alpha)^s}.
\] By Lemma \ref{lemma:2} unless $2m \equiv 0 \bmod w$ the inner
summation is zero and if $2m \equiv 0 \bmod w$ the inner summation is
\[ \sum_{\alpha \in \mathcal V} \frac{(\alpha/\overline{\alpha})^m
\chi_k(\alpha) }{N(\alpha)^s} = w
\frac{L_{K,0}(s,\tau^{(m)}\chi_k)}{\zeta(2s,\chi_k)}. \] Therefore if
$2m \equiv 0 \bmod w$
\[ 
\Upsilon^{(m)}(2s) = \frac{w}2 \sum_{r|2d} \sum_{\ell | 2d/r}
\sum_{k|r\ell} \big(r^s \mu(\ell) \mu(k) \big)
\frac{L_{K,0}(s,\tau^{(m)}\chi_k)}{\zeta(2s,\chi_k)}.
\]

By Lemma \ref{lemma:5}, $\tau^{(m)} \chi_r$ is a Hecke character and
is principal only if $m=0$.  Therefore $L_{K,0}(s,\tau^{(0)}\chi_k)$
has a meromorphic continuation to a neighborhood of $\rp{s} \geq 1$
with only a simple pole at $s=1$ and if $m\neq 0$
$L_{K,0}(s,\tau^{(m)}\chi_k)$ has an analytic continuation to a
neighborhood of $\rp{s}\geq 1$.  Since the summations defining
$\Upsilon^{(m)}(2s) $ are finite, the analogous statements are true
for $\Upsilon^{(m)}.$ The proof of Theorem \ref{thm:1} now follows
using Weyl's criterion.

\bibliography{bibliography}

\begin{thebibliography}{1}

\bibitem{MR2216774}
Enrico Bombieri and Walter Gubler.
\newblock {\em Heights in {D}iophantine geometry}, volume~4 of {\em New
  Mathematical Monographs}.
\newblock Cambridge University Press, Cambridge, 2006.

\bibitem{MR0011302}
Beno Eckmann.
\newblock \"{U}ber monothetische {G}ruppen.
\newblock {\em Comment. Math. Helv.}, 16:249--263, 1944.

\bibitem{MR1544392}
E.~Hecke.
\newblock Eine neue {A}rt von {Z}etafunktionen und ihre {B}eziehungen zur
  {V}erteilung der {P}rimzahlen.
\newblock {\em Math. Z.}, 6(1-2):11--51, 1920.

\bibitem{Jacobson:1989kx}
Nathan Jacobson.
\newblock {\em Basic algebra II}.
\newblock W.H. Freeman, New York, 1980.

\bibitem{MR0160763}
Serge Lang.
\newblock {\em Algebraic numbers}.
\newblock Addison-Wesley Publishing Co., Inc., Reading, Mass.-Palo Alto-London,
  1964.

\bibitem{MR0263823}
Jean-Pierre Serre.
\newblock {\em Abelian {$l$}-adic representations and elliptic curves}.
\newblock McGill University lecture notes written with the collaboration of
  Willem Kuyk and John Labute. W. A. Benjamin, Inc., New York-Amsterdam, 1968.

\bibitem{MR1511862}
Hermann Weyl.
\newblock \"{U}ber die {G}leichverteilung von {Z}ahlen mod. {E}ins.
\newblock {\em Math. Ann.}, 77(3):313--352, 1916.

\end{thebibliography}

\noindent\rule{4cm}{.5pt}
\vspace{.25cm}

\vspace{.25cm}
{\small \noindent {\sc Kathleen L. Petersen} \\
{{ Florida State University} \\
 { Department of Mathematics} \\
 { Tallahassee, Florida} \\
 email: {\tt petersen@math.fsu.edu}}}

\vspace{.5cm} 
{\small \noindent {\sc Christopher D. Sinclair} \\
Department of Mathematics \\
 University of Oregon \\
 Eugene, Oregon \\
email: {\tt csinclai@uoregon.edu}}

\end{document}